\documentclass[preprint]{acmconf} 
\usepackage{amssymb,amsmath,array,delarray,latexsym,epsfig} 

\newlength{\sh}
\newlength{\baker}
\newlength{\greg}
\newlength{\fw}
\newlength{\jmr}
\newlength{\jfc}
\newlength{\bernd}
\newlength{\jones}
\newlength{\mati}
\newlength{\tung}
\newlength{\sil}
\newlength{\koi}
\newlength{\gala}
\newtheorem{lemma}{Lemma}

\newtheorem{dfn}{Definition}

\newtheorem{main}{Theorem} 
\newtheorem{thm}[main]{Theorem}
\newtheorem{cor}{Corollary}
\newtheorem{rem}{Remark}

\newcommand{\tf}{{\tilde{f}}}
\newcommand{\twF}{{\tilde{F}}}

\newcommand{\pspa}{{\mathbf{PSPACE}}}

\newcommand{\am}{{\mathbf{AM}}} 

\newcommand{\np}{{\mathbf{NP}}}

\newcommand{\pp}{\mathbf{P}}
\newcommand{\hn}{\mathbf{HN}}

\newcommand{\eps}{\varepsilon}
\newcommand{\cA}{\mathcal{A}}

\newcommand{\cO}{\mathcal{O}}
\newcommand{\supp}{\mathrm{Supp}}

\newcommand{\thth}{{\underline{\mathrm{th}}}}

\newcommand{\st}{ {\underline{ \mathrm{st} } }  }
\newcommand{\nd}{{\underline{\mathrm{nd}}}}

\newcommand{\Q}{\mathbb{Q}}
\newcommand{\R}{\mathbb{R}}
\newcommand{\C}{\mathbb{C}}
\newcommand{\N}{\mathbb{N}}
\newcommand{\Z}{\mathbb{Z}}

\newcommand{\pert}{\mathrm{Pert}}

\newcommand{\res}{\mathrm{Res}}

\newcommand{\Zn}{\Z^n}

\newcommand{\Rn}{\R^n}

\newcommand{\Cn}{\C^n}

\newcommand{\Csn}{{(\C^*)}^n}
\newcommand{\qed}{$\blacksquare$}

\newcommand{\cI}{\mathcal{I}}

\newcommand{\cR}{\mathcal{R}}
\newcommand{\cP}{\mathcal{P}} 

\newcommand{\cS}{\mathcal{S}}
\newcommand{\cC}{\mathcal{C}}
\newcommand{\cU}{\mathcal{U}}
\newcommand{\bO}{\mathbf{O}}
\newcommand{\vol}{\mathrm{Vol}}

\begin{document}

\title{Computing Complex Dimension Faster and Deterministically\\ (Extended 
Abstract)}  

\author{J.\ Maurice Rojas\thanks{\today {} version. This research was 
partially supported by a Hong Kong CERG grant. }  
\medskip\\
Department of Mathematics\\
City University of Hong Kong\\
83 Tat Chee Avenue\\ 
Kowloon, HONG KONG\\
{\tt mamrojas@math.cityu.edu.hk}\\  
{\tt http://math.cityu.edu.hk/\~{}mamrojas } } 

\date{\today} 

\maketitle

\copyrightspace 

\section{Introduction and Main Results}
\label{sec:intro} 
We give a new complexity bound for calculating the complex dimension of an 
algebraic set. Our algorithm is completely deterministic and approaches the 
best recent randomized complexity bounds. We also present some new, 
significantly sharper quantitative estimates on {\bf rational univariate 
representations (RUR)} of roots of polynomial systems. As a corollary of the 
latter bounds, we considerably improve a recent algorithm of Koiran for 
deciding the emptiness of a hypersurface intersection over $\C$, given 
the truth of the Generalized\footnote{  
The {\bf Riemann Hypothesis (RH)} is an 1859 
conjecture equivalent to a sharp quantitative statement on the 
distribution of primes. GRH can be phrased as a generalization of this 
statement to prime ideals in an arbitrary number field. Further background on 
these RH's can be found in \cite{lago,bs}.} 
Riemann Hypothesis (GRH). 

Let $f_1,\ldots,f_m\!\in\!\C[x_1,\ldots,x_n]$, $\pmb{F}\!:=\!(f_1,\ldots,
f_m)$, and let $Z_F$ be the complex zero set of $F$.  
We will first consider the complexity of computing the complex dimension 
$\dim Z_F$ relative to the BSS model over $\C$. (See \cite{bcss} 
for further background on this computational model.) 

First recall the usual notions of input size: With the Turing model, we 
will assume that any input 
polynomial is given as a sum of monomial terms, with all coefficients {\bf 
and} exponents written in, say, base $2$. The corresponding notion of 
{\bf sparse size} is then simply the total number of bits in all coefficients 
and exponents. For example, the sparse size of $x^D_1+ax^3_1+b$ is $\cO(\log 
D+\log a +\log b)$. The sparse size can be extended to the BSS model 
over $\C$ simply by counting just the total number of bits necessary to write 
down the exponents (thus ignoring the size of the coefficients). 

Curiously, efficient {\bf randomization-free} algorithms for  
computing $\dim Z_F$ are hard to find in the literature. So we present such 
an algorithm, with an explicit complexity bound. 
\begin{main}
\label{main:complex}
Let $D$ the maximum of the total degrees of $f_1,\ldots,f_m$. Also 
let $\bO$ be the origin, and $e_1,\ldots,e_n$ the standard basis vectors, in 
$\Rn$. Normalize $n$-dimensional volume $\vol_n(\cdot)$ so that 
the standard $n$-simplex (with vertices $\bO,e_1,\ldots,e_n$)  
has $n$-volume $1$. Finally, let $k$ be the total number of monomial terms in 
$F$, counting reptitions between distinct $f_i$. Then there is a 
deterministic algorithm which computes $\dim Z_F$ within 
$\cO(n^{2.312}{11}^n k V^{7.376}_F)$ arithmetic operations, where 
$V_F\!:=\!\vol_n(Q_F)$ and $Q_F$ is 
the convex hull of\footnote{i.e., smallest convex set in $\R^3$ containing...} 
the union of $\{\bO,e_1,\ldots,e_n\}$ and the set of all exponent vectors of 
$F$. 
\end{main}
\noindent 
Via a height\footnote{The (absolute logarithmic) {\bf height} of an algebraic 
number can be defined as the sparse size of its minimal polynomial. An 
analogous characterization of this important number-theoretic invariant can 
also be given for any algebraic point in $\Cn$ \cite{sil,gregoheight,cool}. } 
estimate from theorem \ref{main:height} later in 
this section one can also derive a similar bound on the bit complexity of 
dimension computation.  We clarify the benefits of our result over earlier 
bounds in section \ref{sub:relatedcomplex} and give an example in section 
\ref{sub:3by3}. The algorithm for theorem 
\ref{main:complex}, and its correctness proof, are 
stated in section \ref{sub:proofcomplex}. 

There is, however, a fundamentally different approach which, given the truth 
of GRH, places a special case of the above problem in an even better 
complexity class. In particular, let $\hn$ denote the 
problem of deciding whether $Z_F$ empty, given that the coefficients 
of all $f_i$ are integers. Then by a recent result of Koiran \cite{hnam}, 
we know that $\hn\!\in\!\am$, given the truth of GRH. 
Koiran's conditional result gives the smallest complexity class known to 
contain $\hn$. (Without GRH, we only know that $\hn\!\in\!\pspa$ 
\cite{pspace}.) Indeed, 
independent of GRH, while 
it is known that $\np\!\subseteq\!\am\!\subseteq\!\pspa$ \cite{papa}, the 
properness of each inclusion is still an open problem.

The simplest summary of Koiran's algorithm is that it uses 
reduction modulo specially selected primes to decide feasibility over $\C$. 
(His algorithm is unique in this respect since all previous 
algorithms for $\hn$ worked primarily in the ring $\C[x_1,\ldots,x_n]/\langle 
F \rangle$.) The key observation behind Koiran's algorithm is that 
an $F$ infeasible (resp.\ feasible) over $\C$ will have roots in $\Z/p\Z$ 
for only finitely many (resp.\ a positive density of) primes $p$. 
(We give an explicit example of an improved version of Koiran's 
algorithm a bit later in section \ref{sub:3by3}.) 

A refined characterization of the difference between positive and zero 
density, which significantly improves Koiran's algorithm, can be given in 
terms of our framework as follows:  
\begin{thm} 
\label{main:koi}
Following the notation above, assume now that 
$f_1,\ldots,f_m\!\in\!\Z[x_1,\ldots,x_n]$ and let $\sigma(F)$ be the maximum 
of $\log|c|$ as $c$ ranges over the coefficients of all the monomial terms of 
$F$. Then there exist $a_F,A_F\!\in\!\N$, with the following properties: 
\begin{itemize}
\item[(a)]{$F$ infeasible over $\C \Longrightarrow$ the reduction of 
$F$ mod $p$ has a root in $\Z/p\Z$ for at most $a_F$ distinct primes $p$. } 
\item[(b)]{Given the truth of GRH, $F$ feasible over $\C \Longrightarrow$ for 
each $t\!\geq\!4963041$, the sequence 
\mbox{$\{A_Ft^3,\ldots,A_F(t+1)^3-1\}$} contains 
a prime $p$ such that the reduction of $F$ mod $p$ has a root in $\Z/p\Z$. } 
\item[(c)]{We have\\ 
$a_F\!=\!\cO(n^3DV_F(4^nD\log D + \sigma(F) +\log m))$ and\\  
$A_F\!=\!O([\frac{e^n}{\sqrt{n}}V_F(\sigma(F)+m(n\log D+\log m))]^4)$. } 
\end{itemize}
In particular, the bit-sizes of $a_F$ and $A_F$ are both 
\mbox{$\cO(n\log D+\log \sigma(F))$} --- 
sub-quadratic in the sparse size of $F$. Simple explicit 
formulae for $a_F$ and $A_F$ appear in remarks \ref{rem:shebanga} and 
\ref{rem:shebangb} of section \ref{sub:proofcomplex}. 
\end{thm}

The proof of theorem \ref{main:koi} is based in part on a new, highly  
refined version of effective univariate reduction. 
\begin{thm} 
\label{main:height}
Following the notation above, and the assumptions of theorem \ref{main:koi}, 
there exist a univariate polynomial $h_F\!\in\!\Z[u_0]$ and 
a point $u_F\!:=\!(u_1,\ldots,u_n)\!\in\Zn$ with the following 
properties: 
\begin{enumerate} 
\setcounter{enumi}{-1}  
\item{The degree of $h_F$ is $\leq\!V_F$.} 
\item{For any irreducible component $W$ of $Z_F$, there is a 
point $(\zeta_1,\ldots,\zeta_n)\!\in\!W$ such that 
$u_1\zeta_1+\cdots+u_n\zeta_n$ is a root of $h_F$. Conversely, 
if $m\!\leq\!n$, all roots of $h_F$ arise this way. } 
\item{$F$ has only finitely many complex roots $\Longrightarrow$ the splitting 
field of $h_F$ over $\Q$ is exactly the field $\Q[x_i \; \; | \; \; 
(x_1,\ldots,x_n)\!\in\!\Cn \text{ \ is \ a \ root \ of \ } F]$. } 
\item{The coefficients of $h_F$ satisfy 
$\sigma(h_F)\!=\!\cO(\frac{e^n}{\sqrt{n}}V_F(\sigma(F)+m(n\log D
+\log m)))$.}  
\item{$m\!\leq\!n \Longrightarrow$ the deterministic arithmetic complexity 
of computing $u_F$, and all the coefficients of $h_F$, is $\cO(n^{1.312}{11}^n 
V^{7.376}_F)$. } 
\item{We have $\log(1+|u_i|)\!=\!\cO(n^2\log D)$ for all $i$.} 
\end{enumerate} 
\end{thm} 
\noindent 
Note that we have thus obtained the existence of points of 
bounded height on the positive-dimensional part of $Z_F$, 
as well as a bound on the height of any point in the 
zero-dimensional part of $Z_F$. Put more simply, via a slight 
variation of the proof of theorem \ref{main:height}, we obtain 
the following useful bound:  
\begin{thm} 
\label{main:size} 
Following the notation of theorem \ref{main:height}, 
any irreducible component $W$ of $Z_F$ contains a point 
$(\zeta_1,\ldots,\zeta_n)$ such that for all $i$, either 
$x_i\!=\!0$ or \mbox{$|\log|x_i||\!=\!\cO(\frac{e^n}{\sqrt{n}}
V_F(\sigma(F)+m(n\log D+\log m)))$. \qed } 
\end{thm}

Our final main result is a refinement of theorem \ref{main:height} which 
will also prove quite useful. 
\begin{thm} 
\cite{gcp}
\label{main:unired} 
Following the notation of theorem \ref{main:height}, 
one can pick $u_F$ and $h_F$ (still satisfying (0)--(5)) so that there 
exist $a_1,\ldots,a_n\!\in\!\N$ and $h_1,\ldots,h_n\!\in\!\Z[u_0]$ with the 
following properties: 
\begin{enumerate} 
\setcounter{enumi}{5}
\item{ The degrees of $h_1,\ldots,h_n$ are all bounded above 
by $V_F$.} 
\item{ For any root $\theta\!=\!u_1\zeta_1+\cdots+u_n\zeta_n$ of $h_F$, 
$\frac{h_i(\theta)}{a_i}\!=\!\zeta_i$ for all $i$. } 
\item{For all $i$, both $\log a_i$ and 
$\sigma(h_i)$ are bounded above by\\ 
$\cO(\frac{e^n}{\sqrt{n}}V^3_F(\sigma(F)+m(n\log D +\log m)))$.} 
\item{$m\!\leq\!n \Longrightarrow$ the deterministic arithmetic complexity of 
computing all the coefficients of $h_1,\ldots,h_n$ is 
$\cO(n^{2.312}{11}^n V^{7.376}_F )$. } 
\end{enumerate}
\end{thm} 

Explicit formulae for all these asymptotic estimates appear below 
in remarks \ref{rem:size}, \ref{rem:denom}, \ref{rem:shebanga}, 
\ref{rem:shebangb}, and \ref{rem:height}. The 
proofs for all our results (except theorem \ref{main:height}), appear in 
section \ref{sub:proofcomplex}. The proof of theorem \ref{main:height} appears 
in the appendix. However, let us first compare our results to 
earlier work. 

\subsection{Related Results Over $\C$} 
\label{sub:relatedcomplex}  
We point out that we 
have tried to balance generality, sharpness, and ease of proof in our 
bounds. In particular, our bounds fill a lacuna in the literature where 
earlier bounds seemed to sacrifice generality for sharpness, or vice-versa. 

To clarify this trade-off, first note that 
$\cI_F\!\leq\!V_F\!\leq\!D^n$, 
where $\cI_F$ is the number of irreducible components of $Z_F$. 
(The first inequality follows immediately from theorem 
\ref{main:height}, while the 
second follows from the observation that $Q_F$ always lies in a 
copy of the standard $n$-simplex scaled by a factor of $D$.) 
So depending on the shape of $Q_F$, and thus 
somewhat on the sparsity of $F$, one can typically expect $V_F$ 
to be much smaller than $D^n$. For example, the $3\times 3$ 
example from section \ref{sub:3by3} below gives $D^n\!=\!13824$ and 
$V_F\!=\!243$. 
More generally, it is easy to see that the factor of improvement can even 
reach $D^{n-1}$, if not more \cite{front}. 

Our algorithm for computing $\dim Z_F$ gives the first deterministic 
complexity bound which is polynomial in $V_F$. In particular, 
while harder problems were already known to admit $\pspa$ complexity bounds, 
the corresponding complexity bounds were either polynomial (or worse) in 
$D^n$, or stated in terms of a non-uniform computational model.\footnote{ For 
example, some algorithms in the literature are stated in terms 
of {\bf arithmetic networks}, where the construction of the 
underlying network is not included in the complexity estimate. }  
Our algorithm for the computation of $\dim Z_F$ thus gives a significant 
speed-up over earlier work. 

For example, via the work of Chistov and Grigoriev from the early 1980's on 
quantifier elimination over $\C$ \cite{chigo}, it is not hard to derive a 
deterministic complexity 
bound of $\cO((mD)^{n^4})$ for the computation of $\dim Z$. More recently, 
\cite{giustiheintz} gave a randomized complexity bound of 
$m^{\cO(1)}D^{\cO(n)}$. Theorem \ref{main:complex} thus clearly 
improves the former bound. Comparison with the latter bound is 
a bit more difficult since the exponential constants and derandomization 
complexity are not explicit in \cite{giustiheintz}. 

As for faster algorithms, one can seek complexity bounds which are polynomial 
in even smaller quantities. For example, if one has an irreducible algebraic 
variety $V\!\subseteq\!\Cn$ of complex dimension $d$, one can 
define its {\bf affine geometric degree}, $\delta(V)$, to be the number of 
points in $V\cap H$ where $H$ is a generic $(n-d)$-flat.
More generally, we can 
define $\delta(Z_F)$ to be the sum of $\delta(V)$ as $V$ ranges over all 
irreducible components of $Z_F$. It then follows (from theorem 
\ref{main:complex} and a consideration of intersection multiplicities) 
that $\cI_F\!\leq\!\delta(Z_F)\!\leq\!V_F$. 
Similarly, one can attempt to use mixed volumes of several polytopes (instead 
of a single polytope volume) to lower our bounds. 

We have avoided refinements of this nature for the sake of simplicity. 
Another reason it is convenient to have bounds in terms of $V_F$ 
is that the computation of $\delta(Z_F)$ is even more subtle 
than the computation of polytopal $n$-volume. For example, when $n$ is 
fixed, $\vol_n(Q)$ can be computed in polynomial time simply by 
triangulating the polytope $Q$ and adding together the volumes of the 
resulting $n$-simplices 
\cite{volcomplex}. However, merely deciding $\delta(Z_F)\!>\!0$ is 
already $\np$-hard for $(m,n)\!=\!(2,1)$, via a result of 
Plaisted \cite{plaisted}. 
As for varying $n$, computing 
$\delta(Z_F)$ is $\#\pp$-hard, while the computation of polytope volumes is 
$\#\pp$-complete. (The latter result is covered in \cite{volcomplex,kls}, 
while the former result follows immediately from the fact that the computation 
of $\delta(Z_F)$ includes the computation of $V_F$ as a special case.) More 
practically, for any fixed $\eps_1,\eps_2\!>\!0$, there 
is an algorithm which runs in time polynomial in the sparse encoding of $F$ 
(and thus polynomial in $n$) which produces a random variable that is within a 
factor of $1-\eps_1$ of $\vol_n(Q_F)$ with probability $1-\eps_2$ \cite{kls}. 
The analogous result for mixed volume is known only for certain families 
of polytopes \cite{gs00}, and the existence of such a result for 
$\delta(Z_F)$ is still an open problem. 

In any event, we point out that improvements in terms of $\delta(Z_F)$ 
for our bounds are possible, and these will be pursued in a later version 
of this paper. Similarly, the exponents in our complexity bounds 
can be considerably lowered if randomization is allowed. Furthermore, Lecerf 
has recently announced a randomized complexity 
bound for computing $\dim Z_F$ which is polynomial in 
$\max_i\{\delta(Z_{(f_1,\ldots,f_i)})\}$ 
\cite{lecerf}.\footnote{ The paper \cite{lecerf} actually solves the  
harder problem of computing an algebraic description of a non-empty 
set of points in every irreducible component of $Z_F$, and distinguishing 
which component each set belongs to.} However, the complexity of 
derandomizing Lecerf's algorithm is not yet clear. 

As for our result on prime densities (theorem \ref{main:koi}), part (a) 
presents the best current bound polynomial in $V_F$. An earlier 
density bound, polynomial in $D^{n^{\cO(1)}}$ instead, 
appeared in \cite{hnam}. 

Part (b) of theorem \ref{main:koi} appears to be new, and 
makes explicit an allusion of Koiran in \cite{hnam}. 
\begin{rem} 
Pascal Koiran has also given an $\am$ algorithm 
(again depending on GRH) for deciding whether the complex dimension 
of an algebraic set is less than some input constant \cite{koiran}. 
\qed 
\end{rem} 

Regarding our height bound, the only other results stated in polytopal terms 
are an earlier version of theorem \ref{main:height} announced in 
\cite{stoc99}, and independently 
discovered bounds in \cite[Prop.\ 2.11]{cool} and \cite[Cor.\ 
8.2.3]{maillot}. The bound from \cite{cool} applies to a slightly 
different problem, but implies (by intersecting with a generic 
linear subspace with reasonably bounded coefficients)\footnote{ 
Martin Sombra pointed this out in an e-mail to the author.} a bound of  
$\cO((4^nD\log n + n\sigma)V_F)$ for our setting. Their 
bound thus results in a slightly stronger exponential dependence on $n$. The 
bound from \cite[Cor.\ 8.2.3]{maillot} uses Arakelov 
intersection theory, holds only for $m\!=\!n$, and the statement is more 
intricate (involving a sum of several mixed volumes). So it is not yet 
clear when \cite[Cor.\ 8.2.3]{maillot} is better than theorem 
\ref{main:height}. In any case, our result has a considerably simpler 
proof than either of these two alternative bounds: We use only resultants and 
elementary linear algebra and factoring estimates.  

We also point out that the only earlier 
bounds which may be competitive with theorem \ref{main:height}, 
\cite[Prop.\ 2.11]{cool}, and \cite[Cor.\ 8.2.3]{maillot} are polynomial in 
$D^n$ and make various non-degeneracy hypothesis, e.g., 
$m\!=\!n$ and no singularities for $Z_F$ (see \cite{cannyphd} and 
\cite[Thm.\ 5]{gregogap}). As for 
bounds with greater generality, \cite{fgm} gives a height bound for 
general quantifier elimination which, unfortunately, has a factor of the 
form $2^{(n\log D)^{\cO(r)}}$ where $r$ is the number of 
quantifier alternations. 

As for our refinement of theorem \ref{main:height}, 
the approach of RUR for the roots of polynomial systems is not new, and 
even dates back to Kronecker. RUR also goes under the name of 
``effective primitive element theorem'' and important precursors to 
our theorem \ref{main:unired} are stated in \cite{pspace} and 
\cite[Thm.\ 4]{hnam}. Nevertheless, the use of 
{\bf toric resultants} (cf.\ section \ref{sub:proofcomplex}), 
which form the core of our algorithms here, was 
not studied in the context of RUR until the late 1990's (see, e.g., 
\cite{gcp}). Also, theorem \ref{main:unired} appears to be the 
first statement giving bounds on $\sigma(h_i)$ which are polynomial 
in $V_F$. More recently, an algorithm for RUR with randomized complexity 
polynomial in $\max_i\{\delta(Z_{(f_1,\ldots,f_i)})\}$ was derived in 
\cite{gls99}. However, their 
algorithm makes various nondegeneracy assumptions (such as $m\!=\!n$ and that 
$F$ form a complete intersection) and the derandomization 
complexity is not stated. 

As for the factors in our complexity and height estimates which are explicitly 
exponential in $n$ (e.g., $e^n$ and $11^n$), these can be replaced by 
a quantity no worse than $\cO(n^{2.376})$ in certain cases. In 
general, this can be done whenever there exists an expression for a particular 
toric resultant (cf.\ section \ref{sub:proofcomplex}) as a single 
determinant, or the divisor of a determinant, of a matrix of size $\cO(nV_F)$. 
The existence of such formulae has been proven in various cases, e.g., when 
all the Newton polytopes are axis-parallel parallelepipeds \cite{wz}. Also, 
such formulae have been observed (and constructed) experimentally in various 
additional cases of practical interest \cite{emican}. Finding compact formulae 
for resultants is an area of active research which thus has deep implications 
for the complexity of algebraic geometry.  

Finally, we note that 
we have avoided Gr\"obner basis techniques because there are currently 
no known complexity or height bounds polynomial in $V_F$ using 
these methods for the problems we consider. A further complication is that 
there are examples of ideals, generated by polynomials of degree $\leq\!5$ in 
$\cO(n)$ variables, where every Gr\"obner basis has a generator of degree 
$2^{2^n}$ \cite{mm}. This is one obstruction to deriving sharp explicit 
complexity bounds via a naive application of Gr\"obner bases. Nevertheless, 
we point out that Gr\"obner bases are well-suited for other difficult 
algebraic problems, and their complexity is also an area of active research. 

\subsection{A Sparse $\pmb{3\times 3}$ Example}  
\label{sub:3by3}
The solution of sparse polynomial systems is a problem with   
numerous applications outside, as well as inside, mathematics. The 
analysis of chemical reactions \cite{gaterhub} and the computation of 
equilibria in game-theoretic models \cite{mucks} are but two diverse examples. 

More concretely, consider the following system of $3$ polynomial equations in 
$3$ variables: 
\begin{eqnarray}
\label{eq:3by3}
144+2x-3y^2+x^7y^8z^9 &=& 0 \notag \\
-51+5x^2-27z+x^9y^7z^8 &=& 0 \\
7-6x+8x^8y^9z^7-12x^8y^8z^7 &=& 0. \notag 
\end{eqnarray} 
Let us see if the system (\ref{eq:3by3}) has any {\bf complex} 
roots and, if so, count how many there are. 

Note that the total degree\footnote{ The {\bf total degree} of a polynomial is 
just the maximum of the sum of the exponents in any monomial term of the 
polynomial.} of each polynomial above is 24. By 
an 18$^\thth$-century theorem of \'Etienne B\'ezout \cite{shafa}, 
we can bound from above the number of complex roots of (\ref{eq:3by3}), 
assuming 
this number is finite, by $24\cdot 24 \cdot 24 = \mathbf{13824}$. 
However, a more precise 20$^\thth$-century bound can be obtained 
by paying closer attention to the monomial term structure 
of (\ref{eq:3by3}): Considering the 
convex hull of the exponent vectors of each equation in 
(\ref{eq:3by3}), one obtains three tetrahedra. 
These are the {\bf Newton polytopes} of (\ref{eq:3by3}), and 
their {\bf mixed volume}, by a beautiful theorem of 
David N.\ Bernshtein from the 1970's \cite{bernie}, turns out to be a 
much better upper bound on the number of complex roots (assuming there 
are only finitely many). For our polynomial system (1), this bound 
is {\bf 145}. 

Now to decide whether (\ref{eq:3by3}) has any complex roots, 
we can attempt to find a univariate polynomial whose roots 
are some simple function of the roots of (\ref{eq:3by3}). {\bf Elimination 
theory} allows one to do this, and a particularly effective combinatorial 
algorithm is given in theorem \ref{main:complex}. For example, the roots of 
\tiny
$\pmb{P(u)}:= 268435456 u^{145}   -138160373760 u^{137}    
-30953963520 u^{130} +\cdots
-\underline{2947435596503653060289376000} u^{44}+\cdots 
-48803823903916800 u^2  + 8681150210659989300$ 
\normalsize
are exactly those numbers of the form $\alpha\beta\gamma$, where 
$(\alpha,\beta,\gamma)$ ranges over all the roots of (\ref{eq:3by3}) in $\C^3$. 
The above {\bf univariate reduction} thus tells us that our example indeed has 
finitely many complex roots --- exactly 145, in fact. The above 
polynomial took less than $13$ seconds to compute using a 
naive application of {\bf resultants} and factorization on the 
computer algebra system {\tt Maple}. Interestingly, computing 
the same univariate reduction via a naive application 
of {\bf Gr\"obner bases} (on the same machine with the same version of {\tt 
Maple}) takes over $3$ hours and $51$ minutes.

Admittedly, computing polynomials like the one above can be an 
unwieldy approach to deciding whether (1) has a complex root. An alternative 
algorithm, discovered by Pascal Koiran in \cite{hnam} and improved via theorem 
\ref{main:koi}, makes a remarkable 
simplification depending on conjectural properties of the distribution of 
prime ideals in number fields.

For instance, an unoptimized implementation of this alternative algorithm 
would run as follows on our example: 
\begin{itemize} 
\item[{\bf Assumption 1}]{ The truth of the Generalized Riemann Hypothesis 
(GRH). } 
\item[{\bf Step 1}]{ Pick a (uniformly distributed) random integer  
$t\!\in\!\{10^7,\ldots,10^7+2\cdot 10^{11}\}$.} 
\item[{\bf Step 2}]{ Via an oracle in $\np$, decide if there is a prime 
\mbox{$p\!\in\!\{8\cdot 10^{20}\cdot t^3,\ldots,8\cdot 10^{20}\cdot 
(t+1)^3-1\}$} such 
that the mod $p$ reduction of (\ref{eq:3by3}) has a root in $\Z/p\Z$. If so, 
declare  that (\ref{eq:3by3}) has a complex root. Otherwise, declare that 
(\ref{eq:3by3}) has no complex root. \qed } 
\end{itemize} 

The choice of the constants above was assisted via our preceding theorems.  
In particular, the constants 
are simply chosen to be large enough to guarantee that, under GRH, the 
algorithm never fails (resp.\ fails with probability $\leq\!\frac{1}{3}$) 
if (\ref{eq:3by3}) has a complex root (resp.\ does not have a 
complex root). Thus, for our example, the algorithm above will 
always give the right answer regardless of the random choice in 
Step 1.  Note also that while the prime we seek above may be quite large, 
the number of {\bf digits} needed to write any such prime is at most 
${\bf 55}$ --- not much bigger than 53, which is the total number 
of digits needed to write down the coefficients and exponent vectors 
of (\ref{eq:3by3}). 
For the sake of completeness, we observe  
that the number of real (resp.\ rational) roots of (1) is exactly {\bf 11} 
(resp.\ {\bf 0}). 
 
\section{Proofs of Our Results Over $\C$: Theorems \ref{main:complex}, 
\ref{main:height}, \ref{main:size}, \ref{main:unired}, and 
\ref{main:koi}} 
\label{sub:proofcomplex}
While our proof of theorem \ref{main:koi} will not directly 
require knowledge of resultants, our proofs of theorems \ref{main:complex}, 
\ref{main:height}, \ref{main:size}, and \ref{main:unired} are based on the 
{\bf toric resultant}.\footnote{Other 
commonly used prefixes for 
this modern generalization of the classical resultant \cite{vdv} include:
sparse, mixed, sparse mixed, $\cA$-, $(\cA_1,\ldots,\cA_k)$-, and Newton. 
Resultants actually date back to work Cayley and Sylvester in the 
19$^\thth$ century, but the toric resultant incorporates some 
combinatorial advances from the late 20$^\thth$ century. } 
This operator allows us to reduce all the computational algebraic geometry we 
will encounter to matrix and univariate polynomial arithmetic, with almost no 
commutative algebra machinery. 
\begin{rem} 
Another advantage of using toric resultants is their algorithmic uniformity. 
Furthermore, since our algorithms reduce to standard matrix arithmetic in a 
particularly structured way, parallelizing is quite straightforward. \qed 
\end{rem} 
Since we do not have the space to give a full introduction to 
resultants we refer the reader to \cite{emiphd,gkz94,introres} 
for further background. The necessary facts we need are all 
summarized in the appendix of this paper. In what follows, we let 
$[j]\!:=\!\{1,\ldots,j\}$. 

\subsection{The Proof of Theorem \ref{main:complex}}
Our algorithm can be stated briefly as follows: 
\begin{itemize} 
\item[{\bf Step 0}]{If $f_i$ is indentically $0$ for all $i$, 
declare that $Z_F$ has dimension $n$ and stop. Otherwise, 
let $i\!:=\!n-1$. } 
\item[{\bf Step 1}]{For each $j\!\in\![2k+1]$, 
compute an $(i+1)n$-tuple of integers\\ $(\eps_1(j),\ldots,\eps_n(j),
\eps_{(1,1)}(j),\ldots,\eps_{(i,n)}(j))$ via lemma \ref{lemma:probe} and 
the polynomial system (\ref{eq:probe}) below.} 
\item[{\bf Step 2}]{Via theorem \ref{main:height}, 
check if the polynomial system  
\begin{eqnarray} 
\label{eq:probe}  
\eps_1(j)f_1+\cdots+\eps_1(j)^mf_m \mbox{\hspace{6cm}} & & \notag \\
+\eps_1(j)^{m+1}l_1+\cdots+\eps_1(j)^{m+i} l_i = 0 \mbox{\hspace{3cm}} & & 
 \notag \\ 
\vdots \mbox{\hspace{4cm}(2)\hspace{2cm}} & & \notag \\ 
\eps_n(j)f_1+\cdots+\eps_n(j)^mf_m \mbox{\hspace{6cm}} & & \notag \\ 
+\eps_n(j)^{m+1}l_1+\cdots+\eps_n(j)^{m+i} l_i = 0 \mbox{\hspace{3cm}} & & 
\notag 
\end{eqnarray} 
has a root for more than half of the $j\!\in\![2k+1]$, 
where $l_t\!:=\!\eps_{(t,1)}x_1+\cdots+\eps_{(t,n)}x_n$ for 
all $t$. } 
\item[{\bf Step 3}]{If so, declare that $Z_F$ has dimension $i$ and stop. 
Otherwise, if $i\!\geq\!1$, set $i\mapsto i-1$ and go to Step 1.} 
\item[{\bf Step 4}]{ Via theorem \ref{main:unired} and a univariate gcd 
computation, check if the system (\ref{eq:probe}) has a root which is 
also a root of $F$.} 
\item[{\bf Step 5}]{ If so, declare that $Z_F$ has dimension $0$ 
and stop. Otherwise, 
declare $Z_F$ empty and stop.} 
\end{itemize} 
\noindent  
Via the lemma and theorem applied above, we see that Step 2 gives a ``yes''  
answer iff the intersection of $Z_{\twF}$ with a generic 
codimension $i$ flat is finite (and nonempty), where $\twF$ is an 
$n$-tuple of generic 
linear combinations of the $f_i$. Thus Step 2 gives a 
``yes'' answer iff $\dim Z_{\twF}\!=\!i$. 
Lemma \ref{lemma:gh} below tells us that $\dim Z_F\!=\!\dim Z_{\twF}$ if 
$\dim Z_F\!\geq\!1$. Otherwise, Step 5 correctly decides whether $Z_F$ is 
empty whenever $Z_F$ is finite. Thus the algorithm is correct. 

As for the complexity of our algorithm, letting $\cS$ (resp.\ $\cU$, $\cU'$) 
be the complexity bound from lemma \ref{lemma:probe} (resp.\ theorems  
\ref{main:height} and \ref{main:unired}), we immediately obtain a 
deterministic arithmetic complexity bound of\\ 
$\cS+(n-2)(2k+1)\cU+\cU'+kV_F\cO(V_F\log^2 V_F)\!=$
\mbox{$\cO(k\log k+kn^{2.312}e^{2.376n}V^{7.376}_F+n^{2.312}e^{2.376}
V^{7.376}_F)$}
$=\!\cO(n^{2.312}e^{2.376n}kV^{7.376}_F)$. \qed   

\begin{lemma}
\label{lemma:probe}  
Suppose $G(w,v)$ is a formula of the form 
$\exists x_1\!\in\!\C \cdots \exists x_n\!\in\!\C$\\ 
$(g_1(x,w,v)\!=\!0)\wedge \cdots \wedge (g_m(x,w,v)\!=\!0)$,\\
where $g_1,\ldots,g_m\!\in\!\C[x_1,\ldots,x_n,w_1,\ldots,w_k, 
v_1,\ldots,v_r]$. 
Then there is a sequence $v(1),\ldots,v(2k+1)\!\in\!\C^r$ such that 
for all $w\!\in\!\C^k$, the following statement holds: 
$G(w,v(j))$ is true for at least half of the $j\!\in\![2k+1]  
\Longleftrightarrow G(w,v)$ is true for a Zariski-open set of $v\!\in\!\C^r$. 
Furthermore, this sequence can be computed within $\log \sigma + (k+n+r)\log 
D$ arithmetic operations, where $\sigma$ (resp.\ $D$) is the maximum bit-size 
of any  coefficient (resp.\ maximum degree) of any $g_i$. \qed 
\end{lemma} 
\noindent 
The above lemma is actually just a special case of theorem 5.6 of 
\cite{koiran}. 

\subsection{The Proof of Theorem \ref{main:size}}  
Since we only care about the size of $|x_i|$, we can simply 
pick $u_0\!=\!-1$, $u_i\!=1$, all other $u_j\!=\!0$, and apply 
the polynomial $h_F$ from theorem \ref{main:height}. 
(In particular, differing from the proof of theorem \ref{main:height},  
we need not worry if our choice of $(u_1,\ldots,u_n)$ results in two distinct 
$\zeta\!\in\!Z_F$ giving the same value for $\zeta_1u_1+\cdots+\zeta_nu_n$.)  
Thus, by following almost the same proof as assertion (3) of theorem 
\ref{main:height}, we can beat the height bound from theorem 
\ref{main:height} by a summand of $\cO(n^2V_F\log D)$. \qed 
\begin{rem} 
\label{rem:size}
Via theorem \ref{thm:growth} from the appendix (and a classic root 
size estimate of Cauchy \cite{mignotte}), we easily see that the 
asymptotic bound for $|\log|x_i||$ can be replaced by the 
following explicit quantities:  
$\log\left\{\frac{e^{13/6}}{\pi}\sqrt{m_F+1}\cdot 
2^{V_F}4^{m_F}\sqrt{2}^{V_F} \sqrt{\mu}^{m_F}(c+1)^{m_F} \right\}$
\mbox{if $m\!\leq\!n$, or the $\log$ of 
$\frac{e^{13/6}}{\pi}\sqrt{m_F+1}\cdot 
2^{V_F}4^{m_F}\sqrt{2}^{V_F}$}
$\times\sqrt{\mu}^{m_F}(m(mV_F+1)^{m-1}c+1)^{m_F}$
for $m\!>\!n$. \qed 
\end{rem} 

\subsection{The Proof of Theorem \ref{main:unired}}
All portions, save assertion (8), follow immediately from 
\cite[Main Theorem 2.1]{gcp}. To prove assertion (8), we 
will briefly review the computation of $h_1,\ldots,h_n$ 
(which was already detailed at greater length in \cite{gcp}). 
Our height bound will then follow from some elementary  
polynomial and linear algebra bounds. 

In particular, recall the following algorithm for computing 
$h_1,\ldots,h_n$ (the polynomial $\pert$ used below is defined in 
the appendix):  
\begin{itemize} 
\item[{\bf Step 2}]{If $n\!=\!1$, set $h_1(\theta)\!:=\!\theta$ and stop.
Otherwise, for all $i\!\in\![n]$, let $q^-_i(t)$ be the square-free part of
$\pert_A(t,u_1,\ldots,u_{i-1},u_i-1, u_{i+1},\ldots,u_n)$.}
\item[{\bf Step 3}]{Define $q^\star_i(t)$ to be the square-free part of
$\pert_A(t,u_1,\ldots,u_{i-1},u_i+1,u_{i+1},\ldots,u_n)$ for all
$i\!\in\![n]$.}
\item[{\bf Step 4}]{For all $i\!\in\![n]$ and $j\!\in\!\{0,1\}$, let
$r_{i,j}(\theta)$ be the reduction of $\cR_j(q^-_i(t),
q^\star_i((\alpha+1)\theta-\alpha t))$ modulo $h(\theta)$. }
\item[{\bf Step 5}]{For all $i\!\in\![n]$, define
$g_i(\theta)$ to be the reduction of
$-\theta-\frac{r_{i,1}(\theta)}{r_{i,0}(\theta)}$ modulo $h(\theta)$. 
Then define $a_i$ to be the least positive integer so that 
$h_i(t)\!:=\!a_ig_i\!\in\!\Z[t]$. } 
\end{itemize}

Following the notation of the algorithm above, the polynomial 
$\cR_0(f,g)+\cR_1(f,g)t$ is known as the {\bf first subresultant} of $f$ and 
$g$ and can be computed as follows: Letting 
$f(t)\!=\!\alpha_0+\alpha_1t+\cdots+\alpha_{d_1}t^{d_1}$ and
$g(t)\!=\!\beta_0+\beta_1t+\cdots+\beta_{d_2}t^{d_2}$, consider the following
$(d_1+d_2-2)\times (d_1+d_2-1)$ matrix
\begin{tiny}
\[
\begin{bmatrix}
\beta_0 & \cdots & \beta_{d_2} & 0   & \cdots & 0 & 0 \\
0      & \beta_0 & \cdots & \beta_{d_2} & 0 & \cdots & 0\\
\vdots & \ddots & \ddots &  & \ddots & \ddots & \vdots \\
0      & \cdots & 0 & \beta_0 & \cdots & \beta_{d_2} & 0 \\
0      & 0  & \cdots & 0 & \beta_0 & \cdots & \beta_{d_2} \\
\alpha_0 & \cdots & \alpha_{d_1} & 0 & \cdots & 0 & 0 \\
0     & \alpha_0 & \cdots & \alpha_{d_1} & 0 & \cdots & 0\\
\vdots &  \ddots  & \ddots &    &  \ddots  & \ddots  & \vdots   \\
0      & \cdots & 0 & \alpha_0 & \cdots & \alpha_{d_1} & 0 \\
0      & 0  & \cdots & 0 & \alpha_0 & \cdots & \alpha_{d_1}
\end{bmatrix}
\]
\end{tiny}
\hspace{-\sh}with $d_1\!-\!1$ ``$\beta$ rows'' and $d_2\!-\!1$ ``$\alpha$
rows.'' Let $M^1_1$ (resp.\ $M^1_0$) be the submatrix obtained by
deleting the last (resp.\ second to last) column. We then define 
$\cR_i(f,g)\!:=\!\det(M^1_i)$ for $i\!\in\!\{0,1\}$.

Continuing our proof of Theorem \ref{main:unired}, we see that 
we need only bound the coefficient growth of the intermediate 
steps of our preceding algorithm. Thanks to theorem \ref{thm:growth}, 
this is straightforward: First note that 
$\sigma(q^-_i)\!=\!\log((V_F+1)\cdot 2^{V_F})+\sigma(\bar{h}_F)$, 
where $\bar{h}_F$ is the square-free part of $h_F$. (This follows trivially 
from expressing the coefficients of a univariate polynomial $f(t+1)$ 
in terms of the coefficients of $f(t)$.) Via lemma \ref{lemma:mignotte} 
we then see that $\sigma(\bar{h}_F)\!=\!\log(\sqrt{V_F+1}\cdot 2^{V_F})+
\sigma(h_F)$, and thus $\sigma(q^-_i)\!=\!\cO(\sigma(h_F))$. 
Similarly, $\sigma(q^\star_i)\!=\!\cO(\sigma(h_F))$ as well. 

To bound the coefficient growth when we compute $r_{i,j}$ note 
that the coefficient of $t_i$ in\\ 
$q^\star_i(2\theta-t)$ 
is exactly $(-1)^i\sum^d_{j=i} \begin{pmatrix}j \\ i\end{pmatrix} 
(2\theta)^j\alpha_j$, where $\alpha_j$ is the coefficient of 
$t^j$ in $q^\star_i(t)$. Thus, via Hadamard's lemma again, 
we see that 
$|r_{i,j}(\theta)|\!\leq\!\left(\sqrt{V_F+1}\cdot 
e^{\sigma(h_F)}\right)^{V_F-1}\times \\ \left(\sqrt{V_F+1}\cdot 
V_F2^{V_F}(2\theta)^{V_F}
e^{\sigma(h_F)}\right)^{V_F-1}$ for all $i,j$. Since $r_{i,j}$ is 
itself a polynomial in $\theta$ of degree $V_F(V_F-1)$, the 
last inequality then easily implies that 
$\sigma(r_{i,j})\!=\!\cO(V_F\sigma(h_F))$. 

To conclude, note that for any univariate polynomials $f,g\!\in\!\Z[t]$ 
with degree $\leq\!D$, $\sigma(fg)\!=\cO(\sigma(f)+\sigma(g)+\log D)$. 
Via long division it also easily follows that the 
quotient $q$ and remainder $r$ of $f/g$ satisfy $aq,ar\!\in\Z[t]$ 
and $\sigma(aq),\sigma(ar)\!=\!\cO(D(\sigma(f)+\sigma(g)))$, for some 
positive integer $a$ with $\log a\!=\!\cO(\sigma(g))$. 

So by assertion (3) of theorem \ref{main:height} we obtain 
$\log(a_i),\sigma(h_i)\!=\!\cO(V^2_F\sigma(h_F))$, which 
implies our desired bound. \qed  
\begin{rem} 
\label{rem:denom}
An immediately consequence of our proof is that the 
asymptotic bound from assertion (8) can be replaced 
by the following explicit bound: 
$V_F\times\\
\left\{(V_F-1)\left[\log\left(V_F(V_F+1)^{4} {64}^{V_F}\right)+
2\sigma(h_F)\right]+\sigma(h_F) \right\}\\+\sigma(h_F)+\log V_F$. 
\qed
\end{rem}  

\subsection{The Proof of Theorem \ref{main:koi}}
\label{sub:koi}  
\noindent 
{\bf Proofs of Parts (a) and (c):} We first recall the following 
useful effective arithmetic Nullstellensatz of Krick, Pardo, and 
Sombra. 
\begin{thm}
\label{thm:cool} 
Suppose $f_1,\ldots,f_m\!\in\!\Z[x_1,\ldots,x_n]$ and 
$f_1\!=\cdots =\!f_m\!=\!0$ has {\bf no} roots in $\Cn$. 
Then there exist polynomials $g_1,\ldots,g_m\!\in\!\Z[x_1,\ldots,x_n]$ 
and a positive integer $a$ such that $g_1f_1+\cdots +g_mf_m\!=\!a$. 
Furthermore, 
$\log a\!\leq\!2(n+1)^3D V_F[\sigma(F)+\log m + 2^{2n+4}D\log(D+1)]$. \qed 
\end{thm} 
\noindent 
The above theorem is a portion of corollary 3 from \cite{cool}. 

The proof of part (a) is then almost trivial: By assumption, 
theorem \ref{thm:cool} tells us that the mod $p$ reduction of $F$ 
has a root in $\Z/p\Z \Longrightarrow p$ divides $a$. Since 
the number of divisors of an integer $a$ is no more than 
$1+\log a$ (since any prime power other than $2$ is bounded below by 
$e$), we arrive at our desired asymptotic bound on $a_F$. 
So the first half of (c) is proved. \qed 
\begin{rem}
\label{rem:shebanga}
Following the notation of theorem \ref{main:koi}, 
we thus obtain the following explicit bound:  
$a_F\!\leq\!1+2(n+1)^3D V_F[\sigma(F)+\log m + 2^{2n+4}D\log(D+1)]$. 
\qed 
\end{rem}  

\noindent 
{\bf Proofs of Parts (b) and (c):} Recall the following version of the 
discriminant. 
\begin{dfn} 
\label{dfn:disc}
Given any polynomial
$f(x_1)\!=\!\alpha_0+\alpha_1x_1+\cdots+\alpha_Dx^D_1\!\in\!\Z[x_1]$
with all $|\alpha_i|$ bounded above by some integer $c$, define the
{\bf discriminant of} $\mathbf{f}$, $\pmb{\Delta_f}$, to be
$\frac{(-1)^{D(D-1)/2}}{\alpha_D}$ times the following
$(2D-1)\times (2D-1)$ determinant:
\begin{tiny}
\[ \det
\begin{bmatrix}
\alpha_0 & \cdots & \alpha_D & 0 & \cdots & 0 & 0 \\
0     & \alpha_0 & \cdots & \alpha_D & 0 & \cdots & 0\\
\vdots &  \ddots  & \ddots &    &  \ddots  & \ddots  & \vdots   \\
0      & \cdots & 0 & \alpha_0 & \cdots & \alpha_D & 0 \\
0      & 0  & \cdots & 0 & \alpha_0 & \cdots & \alpha_D \\
\alpha_1 & \cdots & D\alpha_D & 0 & \cdots & 0 & 0 \\
0     & \alpha_1 & \cdots & D\alpha_D & 0 & \cdots & 0\\
\vdots &  \ddots  & \ddots &    &  \ddots  & \ddots  & \vdots   \\
0      & \cdots & 0 & \alpha_1 & \cdots & D\alpha_D & 0 \\
0      & 0  & \cdots & 0 & \alpha_1 & \cdots & D\alpha_D 
\end{bmatrix},
\]
\end{tiny} 
\noindent
\mbox{}\hspace{-.15cm}where the first $D-1$ (resp.\ last $D$) rows correspond 
to the coefficients of $f$ (resp.\ the derivative of $f$). \qed
\end{dfn} 

Our proof of part (b) begins with the following observation. 
\begin{thm} 
\label{thm:oyster} 
Suppose $f\!\in\!\Z[x_1]$ is a square-free polynomial of 
degree $D$ with exactly $i_f$ factors over $\Q[x_1]$. Let $N_f(t)$ denote the 
{\bf total} number of distinct roots of the mod $p$ reductions of $f$ in 
$\Z/p\Z$, counted over all primes $p\!\leq\!t$. Then the truth of GRH 
implies that $\left|i_f\pi(t)-N_f(t)\right|\!<\!2\sqrt{t}(D\log t+\log 
|\Delta_f|) +D\log |\Delta_f|$, 
for all $t\!>\!2$. \qed  
\end{thm} 
\noindent 
A slightly less explicit version of the above theorem appeared 
in \cite[Thm.\ 9]{hnam}, and the proof is almost the same as that of an  
earlier result of Adleman and Odlyzko for the case $i_f\!=\!1$ 
\cite[Lemma 3]{amo}. (See also \cite{weinberger}.) The only new 
ingredient is an explicit version of the effective Chebotarev density theorem 
due to Oesterl\'e \cite{oyster}. (Earlier versions of theorem \ref{thm:oyster} 
did not state the asymptotic constants explicitly.) 

The proof of part (b) is essentially a chain of elementary 
analytic bounds which flows from applying theorem \ref{thm:oyster} 
to the polynomial $h_F$ from theorem \ref{main:complex}. However, a 
technicality which must be considered is that $h_F$ might not be 
square-free (i.e., $\Delta_{h_F}$ may vanish). This is easily taken care of by 
an application of the following immediate corollary of lemmata  
\ref{lemma:multi} and \ref{lemma:mignotte}. 
\begin{cor} 
\label{cor:disc}
Following the notation above, let $g$ be the square-free part 
of $f$ and let $D'$ be the degree of $g$. Then 
$\log |\Delta_g|\!\leq\!D'(D\log 2+\log(D'+1)+\log c)$. \qed 
\end{cor}

Another technical lemma we will need regards the existence of sufficiently 
many primes interleaving a simple sequence. 
\begin{lemma} 
\label{lemma:pain} 
The number of primes in the open interval $(At^3,A(t+1)^3)$ 
is at least $\lfloor \frac{1}{12}\cdot\frac{At^2} {\log t+\log A}\rfloor$, 
provided $A,t\!>\!e^5\!\approx\!148.413\ldots$ 
\end{lemma} 
\noindent 
This lemma follows routinely (albeit a bit tediously) from theorem 8.8.4 of 
\cite{bs}, which states that for all $t\!>\!5$, the $t^\thth$ prime lies in 
the open interval $(t\log t,t(\log t+\log\log t))$. 

Our main strategy for proving part (b) is thus the following: 
Let $N_F$ be the obvious analogue of $N_f$ for {\bf systems} of polynomials. 
We will then attempt to find constants $t_0$ and $A_F$ such that 
$N_F(A_F(t+1)^3-1)-N_F(A_Ft^3)\!>\!1$ for all $t\!\geq\!t_0$. 

Via theorems \ref{main:height} and \ref{main:unired}, and a consideration of 
the primes dividing the $a_i$ (the denominators in our rational 
univariate representation of $Z_F$), it immediately follows that 
$|N_F(t)-N_{h_F}(t)|\!\leq\!V_F\sum^n_{i=1}(\log a_i+1)$, for all $t\!>\!0$. 
We are now ready to derive an inequality whose truth will 
imply $N_F(A_F(t+1)^3-1)-N_F(A_Ft^3)\!>\!1$: 
By theorem \ref{thm:oyster}, lemma \ref{lemma:pain}, the triangle inequality, 
and some elementary 
estimates on $\log t$, $t^3$, and their derivatives, it suffices to 
require that $A_Ft^2$ strictly exceed $12(\log A_F+\log t)$ 
times the following quantity: 
$2(1+\sqrt{2})\sqrt{3A_Ft^3}[V_F(\log(3A_Ft^3)+1)+\log |\Delta_g|]+\\
V_F\left(\log|\Delta_g| + \sum^n_{i=1}\log a_i +n\right)+1$, 
for all $t\!>\max\{t_0,e^5\}$, where $g$ denotes the square-free part of 
$h_F$. (Note that we also used the fact that $i_g\!\geq\!1$.) 

A routine but tedious estimation then shows that 
we can actually take $t_0\!=\!1296(\frac{1+\log 3}{3}+\log 1296)\!
\approx\!4963040.506...$, and $A_F$ as in the statement of part (b). 
So part (b) is proved at last. Careful accounting of the estimates then 
easily yields the explicit upper bound for $A_F$ we state below. So the final 
half of part (c) is proved as well. \qed 
\begin{rem} 
\label{rem:shebangb}
The constant $1296(\frac{1+\log 3}{3}+\log 1296)$ arises from trying to find 
the least $t$ for which $t^2\!\geq\!\alpha \log^4t$, where, roughly 
speaking, $\alpha$ ranges over the constants listed in the 
expressions for $A_F,B_F,C_F,D_F$ below:\\ 
\mbox{$A_F\!\leq\!\lceil 1296B^2_F\log^4B_F+36C^2_F\log^2C_F+2D_F\log D_F
\rceil$,} 
where 
$B_F\!:=\!72\sqrt{3}(1+\sqrt{2}) V_F$, 
$C_F\!:=\!24\sqrt{3}(1+\sqrt{2})\log|\Delta_g| +2$, 
and 
$D_F\!:=\!12V_F\left(\log|\Delta_g|+\sum^n_{i=1}\log a_i+n\right)+13$. 
\qed 
\end{rem} 

\section{Acknowledgements} 
The author thanks Felipe Cucker, Steve Smale, and Martin Sombra  
for useful discussions, in person and via e-mail. 

\footnotesize
\bibliographystyle{acm}

\begin{thebibliography}{A}

\bibitem[AO83]{amo} Adleman, Leonard and Odlyzko, Andrew, {\it
``Irreducibility Testing and Factorization of Polynomials,''} 
Mathematics of Computation, 41 (164), pp.\ 699--709, 1983. 

\bibitem[BS96]{bs} Bach, Eric and Shallit, Jeff, {\it
Algorithmic Number Theory, Vol.\ I: Efficient Algorithms,} 
MIT Press, Cambridge, MA, 1996. 

\bibitem[Ber75]{bernie} Bernshtein, D. N., {\it ``The Number of 
Roots of a System of Equations,"} Functional Analysis and its Applications 
(translated from Russian), Vol. 9, No. 2, (1975), 
pp.\ 183--185.

\bibitem[BP94]{binipan} Bini, Dario and Pan, Victor Y. {\it Polynomial and
Matrix Computations, Volume 1: Fundamental Algorithms,} Progress in
Theoretical Computer Science, Birkh\"auser, 1994.

\bibitem[BCSS98]{bcss} Blum, L., Cucker, F., Shub, M., Smale, S., {\it
Complexity and Real Computation,} Springer-Verlag, 1998.

\bibitem[BZ88]{buza} Burago, Yu. D. and Zalgaller, V. A., {\it
Geometric Inequalities,} Grundlehren der mathematischen Wissenschaften 285,
Springer-Verlag (1988).

\bibitem[Can87]{cannyphd} Canny, John F., {\it ``The Complexity of 
Robot Motion Planning Problems,''} ACM Doctoral Dissertation Award 
Series, ACM Press (1987). 

\bibitem[Can88]{pspace} \underline{\hspace{\jfc}}, {\it ``Some Algebraic
and Geometric Computations in PSPACE,''} Proc.\ 20$^\thth$ ACM
Symp.\ Theory of Computing, Chicago (1988), ACM Press.

\bibitem[CG84]{chigo} Chistov, A. L., and Grigoriev, Dima Yu, {\it
``Complexity of Quantifier Elimination in the Theory of Algebraically
Closed Fields,''} Lect.\ Notes Comp.\ Sci.\ 176, Springer-Verlag (1984).

\bibitem[EC93]{emican} Emiris, Ioannis Z.\ and Canny, John, 
{\it ``Efficient Incremental Algorithms for the Sparse Resultant 
and Mixed Volume,''} J.\ Symbolic Comput.\ 20 (1995), no.\ 2, pp.\ 117--149.  

\bibitem[Emi94]{emiphd} Emiris, Ioannis Z., {\it ``Sparse Elimination and
Applications in Kinematics,''} Ph.D. dissertation, Computer Science
Division, U. C. Berkeley (December, 1994), available on-line at {\tt
http://www.inria.fr/saga/emiris}.

\bibitem[EM99]{emimoumat} Emiris, Ioannis Z.\ and Mourrain, Bernard, 
{\it ``Matrices in Elimination Theory,''} J. of Symbolic Computation, 
28(1\&2):3-44, 1999. 

\bibitem[EP99]{emipanmat} Emiris, Ioannis Z.\ and Pan, Victor, 
{\it ``Techniques for Exploiting Structure in Matrix Formulae of the Sparse 
Resultant,''} Toeplitz matrices: structures, algorithms and applications 
(Cortona, 1996), Calcolo 33 (1996), no. 3-4, 353--369 (1998). 

\bibitem[FGM90]{fgm} Fitchas, N., Galligo, A., and Morgenstern, J., 
{\it ``Precise Sequential and Parallel Complexity Bounds for Quantifier 
Elimination Over Algebraically Closed Fields,''} Journal of Pure and 
Applied Algebra, 67:1--14, 1990. 

\bibitem[GH99]{gaterhub} Gatermann, Karin and Huber, Birk, {\it 
``A Family of Sparse Polynomial Systems Arising in Chemical Reaction 
Systems,''} Preprint ZIB (Konrad-Zuse-Zentrum f\"ur Informationstechnik 
Berlin) SC-99 27, 1999. 

\bibitem[GKZ94]{gkz94} Gel'fand, I. M., Kapranov, M. M., and 
Zelevinsky, A. V., {\it Discriminants, Resultants and Multidimensional 
Determinants,} Birkh\"auser, Boston, 1994. 

\bibitem[GH93]{giustiheintz} Giusti, Marc and Heintz, Joos, 
{\it ``La d\'etermination des points isol\'es et la dimension 
d'une vari\'et\'e alg\'ebrique peut se faire en temps polynomial,''} 
Computational Algebraic Geometry and Commutative Algebra (Cortona, 
1991), Sympos.\ Math.\ XXXIV, pp.\ 216--256, Cambridge University 
Press, 1993. 

\bibitem[GLS99]{gls99} Giusti, M., Lecerf, G., and Salvy, B.,
{\it ``A Gr\"obner-Free Alternative to Polynomial System Solving,''}
preprint, TERA, 1999.

\bibitem[GK94]{volcomplex} Gritzmann, Peter and Klee, Victor,
{\it ``On the Complexity of Some Basic Problems in Computational
Convexity II: Volume and Mixed Volumes,''} Polytopes: Abstract,
Convex, and Computational (Scarborough, ON, 1993), pp.\
373--466, NATO Adv.\ Sci.\ Inst.\ Ser.\ C Math.\ Phys.\ Sci.,
440, Kluwer Acad.\ Publ., Dordrecht, 1994.

\bibitem[GS00]{gs00} Gurvits, Leonid and Samorodnitsky, Alex, {``A 
Deterministic Polynomial-Time Algorithm for Approximating Mixed 
Discriminant and Mixed Volume,''} Proceedings of STOC 2000, ACM Press, 2000. 

\bibitem[KLS97]{kls} Kannan, R., Lovasz, L, and Simonovitz, M., 
{\it ``Random Walks and an $\cO^*(n^5)$ Volume Algorithm 
for Convex Bodies,"} 
Random Structures Algorithms, {\bf 11} (1997), no.\ 1, pp.\ 1--50.

\bibitem[Koi96]{hnam} Koiran, Pascal, {\it ``Hilbert's Nullstellensatz 
is in the Polynomial Hierarchy,''} DIMACS Technical Report 96-27, 
July 1996. ({\bf Note:} This preprint considerably improves the published 
version which appeared in Journal of Complexity in 1996.)  

\bibitem[Koi97]{koiran} \underline{\hspace{\koi}}, {\it ``Randomized and 
Deterministic Algorithms for the Dimension of Algebraic Varieties,''} 
Proceedings of the 38$^\thth$ Annual IEEE Computer Society 
Conference on Foundations of Computer Science (FOCS), 
Oct.\ 20--22, 1997, ACM Press. 

\bibitem[KPS00]{cool} Krick, T., Pardo, L.-M., and Sombra, M., {\it 
``Sharp Arithmetic Nullstellensatz,''} submitted for publication, 
also downloadable from {\tt http://xxx.lanl.gov/abs/math.AG/9911094}. 

\bibitem[LO77]{lago} Lagarias, Jeff and Odlyzko, Andrew,
{\it ``Effective Versions of the Chebotarev Density Theorem,''}
Algebraic Number Fields: $L$-functions and Galois Properties
(Proc.\ Sympos.\ Univ.\ Durham, Durham, 1975),  409--464,
Academic Press, London, 1977.

\bibitem[Lec00]{lecerf} Lecerf, Gr\'egoire, {\it ``Computing an 
Equidimensional Decomposition of an Algebraic Variety by Means 
of Geometric Resolutions,''} submitted to the proceedings of the 
International Symposium on Symbolic Algebra and Computation 
(ISSAC) 2000. 

\bibitem[Mai00]{maillot} Maillot, Vincent, {\it 
``G\'eom\'etrie D'Arakelov Des Vari\'et\'es Toriques 
et Fibr\'es en Droites Int\'egrables,''} M\'em.\ Soc.\ Math.\ France, 
to appear. 

\bibitem[Mal00a]{gregogap} Malajovich-Mu\~noz, Gregorio, {\it 
``Condition Number Bounds for Problems with Integer Coefficients,''} 
Journal of Complexity, to appear. 

\bibitem[Mal00b]{gregoheight} \underline{\hspace{\greg}}, {\it ``Transfer 
Theorems for the $\pp\!\neq\!\np$ Conjecture,''} Journal of Complexity, 
to appear.

\bibitem[MM82]{mm} Mayr, E. and Meyer, A., {\it ``The 
Complexity of the Word Problem for Commutative Semigroups 
and Polynomial Ideals,''} Adv.\ Math.\ {\bf 46}, 305--329, 1982. 

\bibitem[MM95]{mucks} McKelvey, Richard D., and McLennan, Andrew,
{\it ``The Maximal Number of Regular Totally Mixed Nash Equilibria,''}
preprint, Department of Economics, University of Minnesota, 1995.    

\bibitem[Mig92]{mignotte} Mignotte, Maurice, {\it
Mathematics for Computer Algebra,} translated from the 
French by Catherine Mignotte, Springer-Verlag, New York, 1992.  

\bibitem[Oes79]{oyster} Oesterl\'e, Joseph, {\it ``Versions Effectives 
du Th\'eor\`eme de Chebotarev sous l'Hypoth\`ese de 
Riemann G\'en\'eralis\'ee,''} 
Ast\'erisque {\bf 61} (1979), pp.\ 165--167. 

\bibitem[Pap95]{papa} Papadimitriou, Christos H., {\it Computational 
Complexity,} Addison-Wesley, 1995. 

\bibitem[Pla84]{plaisted} Plaisted, David A., {\it ``New NP-Hard and 
NP-Complete Polynomial and Integer Divisibility Problems,''} 
Theoret.\ Comput.\ Sci.\ 31 (1984), no.\ 1--2, 125--138. 

\bibitem[Roj99a]{jpaa} Rojas, J.\ Maurice, {\it ``Toric
Intersection Theory for Affine Root Counting,''} Journal of Pure and
Applied Algebra, vol.\ 136, no.\ 1, March, 1999, pp.\ 67-100.

\bibitem[Roj99b]{stoc99} \underline{\hspace{\jmr}}, {\it ``On the Complexity 
of Diophantine Geometry in Low Dimensions,''} Proceedings of the 
31$^\st$ Annual ACM Symposium on Theory of Computing (STOC '99, May 1-4, 1999,
Atlanta, Georgia), 527-536, ACM Press, 1999. 

\bibitem[Roj99c]{gcp} \underline{\hspace{\jmr}}, {\it ``Solving Degenerate 
Sparse Polynomial Systems Faster,''} Journal of Symbolic Computation, 
vol.\ 28 (special issue on
elimination theory), no.\ 1/2, July and August 1999, pp.\ 155--186.

\bibitem[Roj00a]{front} \underline{\hspace{\jmr}}, {\it ``Low-Dimensional
Varieties and the Frontier to Tractability,''} Contemporary Mathematics,
Proceedings of a Conference on Hilbert's Tenth Problem and Related Subjects
(University of Gent, November 1-5, 1999), edited by
Jan Denef, Leonard Lipschitz, Thanases Pheidas, and Jan Van
Geel, AMS Press.

\bibitem[Sch80]{schwartz} Schwartz, J., {\it ``Fast Probabilistic 
Algorithms for Verification of Polynomial Identities,''} 
J.\ of the ACM 27, 701--717, 1980.

\bibitem[Sha94]{shafa} Shafarevich, Igor R., {\it Basic 
Algebraic Geometry I,} second edition, Springer-Verlag (1994).

\bibitem[Sil95]{sil} Silverman, Joseph H., {\it
The Arithmetic of Elliptic Curves,} corrected
reprint of the 1986 original, Graduate
Texts in Mathematics 106, Springer-Verlag (1995).

\bibitem[Stu94]{combiresult} Sturmfels, Bernd, {\it ``On the Newton
Polytope of the Resultant,''} Journal of Algebraic Combinatorics, 3: 207--236,
1994.

\bibitem[Stu98]{introres} \underline{\hspace{\bernd}}, {\it ``Introduction to
Resultants,''} Applications of Computational Algebraic Geometry 
(San Diego, CA, 1997), 25--39, Proc.\ Sympos.\ 
Appl.\ Math., 53, Amer.\ Math.\ Soc., Providence, RI, 1998.

\bibitem[Van50]{vdv} van der Waerden, B. L., {\it Modern Algebra,}
2$^\nd$ edition, F.\ Ungar, New York, 1950.\footnote{Shamefully, 
the sections on resultants were removed from later editions of this book.} 

\bibitem[WZ94]{wz} Weiman, Jerzy and Zelevinsky, Andrei, {\it ``Multigraded 
Formulae for Multigraded Resultants,''} J. Algebraic Geom. 3 (1994), no. 4, 
pp.\ 569--597. 

\bibitem[Wei84]{weinberger} Weinberger, Peter, {\it ``Finding the Number 
of Factors of a Polynomial,''} Journal of Algorithms, 5:180--186, 1984.

\end{thebibliography}

\normalsize

\section{Appendix: Background on Toric Resultants and The Proof of Theorem 
\ref{main:height}} 
\label{sub:height} 
\noindent
Recall that the {\bf support}, $\pmb{\supp(f)}$, of a polynomial 
$f\!\in\!\C[x_1,\ldots,x_n]$ is simply the set of exponent vectors of the 
monomial terms appearing\footnote{We of course fix an ordering on the 
coordinates of the exponents which is compatible with the usual ordering 
of $x_1,\ldots,x_n$.} in $f$. The support of the 
{\bf polynomial system} $F\!=\!(f_1,\ldots,f_m)$ is simply the $m$-tuple 
$\pmb{\supp(F)}\!:=\!(\supp(f_1),\ldots,\supp(f_m))$. Let 
$\bar{\cA}\!=\!(\cA_1,\ldots,\cA_{m+1})$ be any $(m+1)$-tuple of non-empty 
finite subsets of $\Zn$ and set $\cA\!:=\!(\cA_1,\ldots,\cA_m)$. If we say 
that $F$ has {\bf support contained in} $\cA$ then we simply mean that 
$\supp(f_i)\!\subseteq\!\cA_i$ for all $i\!\in\![m]$. 

\begin{dfn}
\label{dfn:res}  
Following the preceding notation, suppose we can find line 
segments $[v_1,w_1],\ldots,[v_{m+1},w_{m+1}]$  
with $\{v_i,w_i\}\!\subseteq\!\cA_i$ for all $i$ and $\vol_m(L)\!>\!0$, where  
$L$ is the convex hull of $\{\bO,w_1-v_1,\ldots,w_{m+1}-v_{m+1}\}$. Then we 
can associate to $\bar{\cA}$ a unique (up to sign) irreducible polynomial 
$\res_{\bar{\cA}}\!\in\!\Z[c_{i,a} \; | \; i\!\in\![m+1] 
 \ , \ a\!\in\!\cA_i]$ with the following property: If we identify 
$\bar{\cC}\!:=\!(c_{i,a} \; | \; i\!\in\![m+1] \ , \ a\!\in\!\cA_i)$ 
with the vector of coefficients of a polynomial system $\bar{F}$ with support 
contained in $\bar{\cA}$ (and constant coefficients), then $\bar{F}$ has a 
root in $\Csn \Longrightarrow \res_{\bar{\cA}}(\bar{\cC})\!=\!0$. Furthermore, 
for all $i$, the degree of $\res_{\bar{\cA}}$ with respect to the coefficients 
of $f_i$ is no greater than $V_F$. \qed 
\end{dfn} 
\noindent 
That the toric resultant can actually be defined as above is covered 
in detail in \cite{combiresult,gkz94}. 

Another operator much closer to our purposes is the {\bf toric 
perturbation} of $F$. 
\begin{dfn}
\label{dfn:pert}  
Following the notation of definition \ref{dfn:res}, assume further that 
$m\!=\!n$, $\supp(F)\!=\!\cA$, and $\supp(F^*)\!\subseteq\!\cA$. 
We then define $\pert_{(F^*,\cA_{n+1})}(u)\!\in\!\C[u_a \; | \; 
a\!\in\!\cA_{n+1}]$ to be the coefficient of the term of\\  
$\res_{\bar{\cA}}(f_1-sf^*_1,\ldots,f_n-sf^*_n,\sum_{a\in 
\cA_{n+1}}u_ax_a)\\
\in\C[s][u_a \; | 
\; a\!\in\!\cA_{n+1}]$ of {\bf lowest} degree in $s$. \qed 
\end{dfn} 
\noindent 
The geometric significance of $\pert$ can be summarized as follows: 
For a suitable choice of $F^*$, $\cA_{n+1}$, and $\{u_a\}$, 
$\pert$ satisfies all the properties of the polynomial $h_F$ from theorem 
\ref{main:height} in the special case $m\!=\!n$. In essence, $\pert$ is an 
algebraic deformation which allows us to replace the positive-dimensional 
part of $Z_F$ by a finite subset which is much easier to handle.
 
To prove theorems \ref{main:complex}, \ref{main:height}, and \ref{main:unired} 
we will thus need a good complexity  estimate for computing $\res$ and $\pert$. 
\begin{lemma}
\label{lemma:respert} 
Following the notation above, 
let $\cR_F$ (resp.\ $\cP_F$) be the number of deterministic arithmetic 
operations needed to evaluate 
$\res_{\bar{\cA}}$ (resp.\ $\pert_{(F^*,\cA_{n+1})}$) at any point in 
$\C^{k+n+1}$ (resp.\ $\C^{2k+n+1}$), where $\cA\!\subseteq\!\supp(F)$ and 
$\cA_{n+1}\!:=\!\{\bO,
e_1,\ldots,e_n\}$. Also let $r_F$ be the total degree of 
$\res_{\bar{\cA}}$ as a polynomial in the coefficients of $\bar{F}$ and set 
$m_F\!:=\!e^{1/8}\frac{e^n}{\sqrt{n+1}}V_F$. (Note that 
$e^{1/8}\!\leq\!1.3315$.) 
Then $r_F\!\leq\!(n+1)V_F$, 
\mbox{$\cR_F\!\leq\!(n+1)r_F\cO(m^{2.376}_F)\!=\!\cO(n^{0.812}e^{2.376n}
V^{3.376}_F)$,} 
and $\cP_F\!\leq\!(r_F+1)\cR_F+O(r_F\log r_F)\!=\!  
\cO(n^{1.812}e^{2.376n}V^{4.376}_F)$. Furthermore, 
$k\!\leq\!mV_F$. \qed 
\end{lemma}
\noindent 
The bound on $\cR_F$ (resp.\ $\cP_F$) follows directly from \cite{emican} 
(resp.\ \cite{gcp}), as well as a basic complexity result on the 
{\bf inverse discrete Fourier transform} \cite[pg.\ 12]{binipan}. 
The very last bound follows from a simple lattice point count. 

Admittedly, such complexity estimates seem rather mysterious without  
any knowledge of how $\res$ and $\pert$ are computed. So let us 
now give a brief summary: 
The key fact to observe is that, in the best circumstances, one can express 
$\res$ as the determinant of a sparse structured matrix $M_{\bar{\cA}}$  
(a {\bf toric resultant matrix}) whose entries are either $0$ or polynomials 
in the coefficients of $\bar{F}$ \cite{emican,emiphd,emipanmat,emimoumat}. 
In fact, the constant 
$m_F$  in our theorem above is nothing more than an upper bound, 
easily derived from \cite{emican} and \cite{gcp}, on the 
number of rows and columns of $M_{\bar{\cA}}$.  

However, it is more frequent that $\res$ is but a {\bf divisor} of such a 
determinant, and further work must be done. Fortunately, in 
\cite{emican,emiphd}, there are general randomized and deterministic 
algorithms for extracting $\res$. 

\subsection*{The Proof of Theorem \ref{main:height}} 
Curiously, precise estimates on coefficient growth in toric resultants 
are absent from the literature. So we supply such an estimate below. 
In what follows, we use $u_i$ in place of $u_{e_i}$. 
\begin{thm} 
\label{thm:growth} 
Following the notation of lemma \ref{lemma:respert}, suppose the coefficients 
of $F$ (resp.\ $F^*$) have absolute value bounded above by $c$ (resp.\ 
$c^*$) for all $i\!\in\![n]$ and $u_1,\ldots,u_n\!\in\!\C$. 
Also let $\|u\|\!:=\!\sqrt{u^2_1+\cdots+u^2_n}$ and let $\mu$ 
denote the maximal number of monomial terms in any $f_i$. 
Then the coefficient of $u^i_0$ in $\pert_{(F^*,\cA_{n+1})}$ has absolute 
value bounded above by $\frac{e^{13/12}}{\sqrt{\pi}}\sqrt{m_F+1}\cdot 
4^{m_F-i/2}\|u\|^{V_F-i}
(\sqrt{\mu}(c+c^*))^{m_F} 
\begin{pmatrix} V_F\\ i \end{pmatrix}$,
assuming that $\det M_{\bar{\cA}}\!\neq\!0$ under the substitution 
$(F-sF^*,u_0+u_1x_1+\cdots+u_nx_n) \mapsto \bar{F}$.  
(Note also that $\frac{e^{13/12}}{\sqrt{\pi}}\!\leq\!1.66691$.)   
\end{thm} 
\noindent 
{\bf Proof:} Let $c_{ij}$ denote the coefficient of 
$u^i_0s^j$ in $\det M_{\bar{\cA}}$, under the substitution
$(F-sF^*,u_0+u_1x_1+\cdots+u_nx_n) \mapsto \bar{F}$. 
Our proof will consist of computing an upper bound on $|c_{ij}|$, so 
we can conclude simply by maximizing over $j$ and then invoking a 
quantitative lemma on factoring. 

To do this, we first observe that one can always construct 
a toric resultant matrix with exactly $n_F$ rows corresponding 
to $f_{n+1}$ (where $\delta(Z_F)\!\leq\!n_F\!\leq\!V_F$), and the 
remaining rows corresponding to $f_1,\ldots,f_n$. 
(This follows from the algorithms we have already invoked 
in lemma \ref{lemma:respert}.) Enumerating how appropriate collections rows 
and columns can contain $i$ entries of $u_0$ (and $j$ entries 
involving $s$), it is easily verified that $c_{ij}$ is a sum of no more than 
$\begin{pmatrix} V_F \\ i \end{pmatrix} 
\begin{pmatrix} m_F-i \\ j \end{pmatrix}$ 
subdeterminants of $M_{\bar{A}}$ of size no greater than $m_F-i-j$. 
The coefficient $c_{ij}$ also receives similar contributions 
from some larger subdeterminants since the rows of $M_{\bar{\cA}}$ 
corresponding to $f_1,\ldots,f_n$ involve terms of the 
form $\gamma+\eps s$. 

Via lemma \ref{lemma:multi} below, we can then derive 
an upper bound of $\begin{pmatrix} V_F \\ i \end{pmatrix}\begin{pmatrix}
 m_F-i\\ j\end{pmatrix} \|u\|^{V_F-i}\\
(\sqrt{\mu}(c+c^*))^{m_F-j}$ on $|c_{ij}|$. However, what we really need is an 
estimate on 
the coefficient $c_i$ of $u^i_0$ of $\pert_{(F^*,\cA_{n+1})}$, assuming 
the non-vanishing of $\det M_{\bar{\cA}}$. To 
estimate $c_i$, we simply apply lemma \ref{lemma:mignotte} below 
(observing that $\pert_{(F^*,\cA_{n+1})}$ is a divisor 
of an $m_F\times m_F$ determinant) to obtain an upper bound of 
$\sqrt{m_F+1}\times$ 
\mbox{$2^{m_F}\begin{pmatrix} V_F \\ i \end{pmatrix}\max_j\left\{
\begin{pmatrix} m_F-i\\ j\end{pmatrix}\right\} 
\|u\|^{V_F-i}(\sqrt{\mu}(c+c^*))^{m_F}$} 
on $|c_i|$. We can then finish via the elementary inequality $\begin{pmatrix} 
m_F-i \\ j \end{pmatrix}\!\leq\!\frac{e^{13/12}}{\sqrt{\pi}}2^{m_F-i}$,  
valid for all $j$ (which in turn is a simple corollary of Stirling's 
formula). \qed 

A simple result on the determinants of certain symbolic matrices, 
used above, is the following. 
\begin{lemma}
\label{lemma:multi} 
Suppose $A$ and $B$ are complex $N\times N$ matrices, where 
$B$ has at most $N'$ nonzero rows. Then the coefficient of 
$s^j$ in $\det(A+sB)$ has absolute value no greater than 
$\begin{pmatrix} N' \\ j \end{pmatrix}v^{N-j}(v+w)^j$, 
where $v$ (resp.\ $w$) is any upper bound on the 
Hermitian norms of the rows of $A$ (resp.\ $B$). \qed 
\end{lemma} 
\noindent 
The lemma follows easily by reducing to the case $j\!=\!0$, via 
the multilinearity of the determinant. The case $j\!=\!0$ is 
then nothing more than the classical {\bf Hadamard's lemma} 
\cite{mignotte}. 

The lemma on factorization we quoted above is the following. 
\begin{lemma} 
\cite{mignotte}
\label{lemma:mignotte}
Suppose $f\!\in\!\Z[x_1,\ldots,x_N]$ has total degree $D$ and
coefficients of absolute value $\leq\!c$. 
Then $g\!\in\Z[x_1,\ldots,x_N]$ divides $f \Longrightarrow$ the 
coefficients of $g$ have absolute value $\leq\!\sqrt{D+1}\cdot 2^Dc$. \qed
\end{lemma} 

We are now ready to prove theorem \ref{main:height}:\\
{\bf Proof of Theorem \ref{main:height}}:\\
{\bf (The Case $\pmb{m\!=\!n}$):} 
The existence of an $h_F$ satisfying (0)--(5) will follow from 
setting $h_F(u_0)\!:=\!\pert_{(F^*,\cA_{n+1})}(u_0)$ for 
$\cA_{n+1}$ as in lemma \ref{lemma:respert}, $F^*$ as in lemma 
\ref{lemma:fill} below, and picking several $(u_1,\ldots,u_n)$ until a good 
one is found. Assertion (0) of theorem \ref{main:height} thus follows 
trivially. That the conclusion of lemma \ref{lemma:fill} implies assertion 
(1) is a consequence of \cite[Def.\ 2.2 and Main Theorem 2.1]{gcp}.  
 
To prove assertions (1)--(5) together we will then need to pick 
$(u_1,\ldots,u_n)$ subject to a final technical condition. In particular, 
consider the following method: 
Pick $\eps\!\in\![1+\begin{pmatrix} V_F \\ 2\end{pmatrix}]$ and set 
$u_i\!:=\!\eps^i$ for all $i\!\in\![n]$.
The worst that can happen is that  
a root of $h_F$ is the image two distinct points 
in $Z_F$ under the map $(\zeta_1,\ldots,\zeta_n) \mapsto 
u_1\zeta_1+\cdots+u_n\zeta_n$, thus obstructing assertion (2). (Whether this 
happens can easily be checked within $\cO(V_F\log V_F)$ arithmetic 
operations via a gcd calculation detailed in \cite[Sec.\ 5.2]{gcp}, 
after first finding the coefficients of $h_F$.) 
Otherwise, it easily follows from Main Theorems 2.1 and 2.4 of \cite{gcp}  
(and theorem \ref{main:unired} above and theorem \ref{thm:growth} below) 
that $h_F$ satisfies assertions (1)--(3) and (5). 

Since there are at most $\begin{pmatrix} 
V_F \\ 2\end{pmatrix}$ pairs of points $(\zeta_1,\zeta_2)$, 
picking $(u_1,\ldots,u_n)$ as specified above {\bf will} eventually 
give us a good $(u_1,\ldots,u_n)$. The overall arithmetic complexity of our 
search for $u_F$ and $h_F$ is, thanks to lemma \ref{lemma:respert},\\ 
$\left(\begin{pmatrix}V_F \\ 2 \end{pmatrix}+1\right) \cdot 
(V_F\cP_F+\cO(V_F\log V_F))$. This proves assertion (4), and we are done. 
\qed 
\begin{rem} 
Note that we never actually had to compute $V_F$ above: To pick a 
suitable $u$, we simply keep pick choices (in lexicographic order) with 
successively larger and larger coodinates until we find a suitable $u$. \qed 
\end{rem} 

\noindent 
{\bf (The Case $\pmb{m\!<\!n}$):} Take $f_{n+1}\!=\cdots =\!f_m\!=\!f_n$. 
Then we are back in the case $m\!=\!n$ and we are done. \qed 

\noindent 
{\bf (The Case $\pmb{m\!>\!n}$):} Here we employ an old trick: We substitute 
generic linear combinations of $f_1,\ldots,f_m$ for $f_1,\ldots,f_n$. 
In particular, set $\tf_i\!:=f_1+\eps_if_2+\cdots+\eps^{m-1}_if_m$ for 
all $i\!\in\![n]$. It then follows from lemma \ref{lemma:gh} below 
that, for generic $(\eps_1,\ldots,\eps_n)$, $Z_{\twF}$ is the union of $Z_F$ 
and a (possibly empty) finite set of 
points. So by the $m\!=\!n$ case, and taking into account the larger 
value for $c$ in our application of theorem \ref{thm:growth}, we are done. \qed 
\begin{rem} 
\label{rem:height}
Via theorem \ref{thm:growth}, we thus see that the 
asymptotic bound of assertion (3) can be replaced by 
essentially the same explicit quantities as detailed 
in remark \ref{rem:size}. The only difference is 
that we replace $\sqrt{2}$ by  
$\sqrt{n}\left(\begin{pmatrix}V_F \\ 2\end{pmatrix}
+1\right)^n$. This accounts for the slightly larger 
asymptotic estimate. \qed 
\end{rem} 
\begin{lemma} 
\label{lemma:fill} 
Following the notation above 
let $\cA^*_i\!=\!\{\bO,e_1,\ldots,e_n\}\cup\bigcup^n_{j=1}\cA_j$ for all 
$i\!\in\![n]$ and $k^*\!:=\!n\#\cA_1$, where 
$\#$ denotes set cardinality. 
Also let $\cC^*$ be the coefficient vector of $F^*$. 
Then there is an $F^*$ such that (i) 
$\supp(F^*)\!\subseteq\!\cA^*$, (ii) 
$\cC^*\!=\!(1,\ldots,1)$, (iii) $F^*$ has exactly 
$V_F$ roots in $\Csn$ counting multiplicities, and 
(iv) $\det M_{\bar{\cA}}\!\neq\!0$ under the substitution 
$(F-sF^*,u_0+u_1x_1+\cdots+u_nx_n) \mapsto \bar{F}$. \qed 
\end{lemma} 
\noindent 
The above lemma is a paraphrase of \cite[Definition 2.3 and Main Theorem 
2.3]{gcp}. 

\begin{lemma}
\label{lemma:gh}
Following the notation above, let $S\!\subset\!\C$ be any finite set 
of cardinality $\geq\!mV_F+1$. Then there is an 
$(\eps_1,\ldots,\eps_n)\!\in\!S^n$ such that every irreducible component of 
$Z_{\twF}$ is either an irreducible component of $Z_F$ or a point. \qed  
\end{lemma} 
\noindent 
The proof is essentially the same as the first theorem of \cite[Sec.\ 
3.4.1]{giustiheintz}, save that we use part (0) of theorem \ref{main:height} 
in place of B\'ezout's Theorem. 

\end{document}